\documentclass[reqno,oneside,11pt]{amsart}

\usepackage[T1]{fontenc}
\usepackage{times,mathptm}
\usepackage{amssymb,epsfig,verbatim,mathtools,float}
\usepackage[all,cmtip]{xy}
\theoremstyle{plain}

\newtheorem{thm}{Theorem}[section]

\newtheorem{pro}[thm]{Proposition}
\newtheorem{lem}[thm]{Lemma}
\newtheorem{proposition-principale}[thm]{Proposition principale}
\newtheorem{theoalph}{Theorem}

\theoremstyle{definition}

\newtheorem{rem}[thm]{Remark}

\addtocounter{section}{0}              
\numberwithin{equation}{section}       
\begin{document}
\setlength{\baselineskip}{0.54cm}         
\title{On regularizable birational maps}
\date{}
\author{Julie D\'eserti}\thanks{The 
 author was partially supported by the ANR grant Fatou ANR-17-CE40-
0002-01 and the ANR grant Foliage ANR-16-CE40-0008-01.}
\address{Universit\'e C\^ote d'Azur, Laboratoire J.-A. Dieudonn\'e, UMR 7351, Nice, France}
\email{deserti@math.cnrs.fr}



\begin{abstract} 
Bedford asked if there exists a 
birational self map $f$ of the complex 
projective plane such that for any 
automorphism $A$ of the complex projective
plane $A\circ f$ is not conjugate to 
an automorphism. 
In this article we give such a $f$ 
of degree $5$.
\end{abstract}

\maketitle

\section{Introduction}

Denote by $\mathrm{Bir}(\mathbb{P}^k_\mathbb{C})$ 
the group of all birational self maps of
$\mathbb{P}^k_\mathbb{C}$, also called the 
\textsl{$k$-dimensional Cremona group}. 
Let $\mathrm{Bir}_d(\mathbb{P}^k_\mathbb{C})$ be the
algebraic variety of all birational self maps of
$\mathbb{P}^k_\mathbb{C}$ of degree~$d$. 
When $k=2$ and $d\geq 2$ these varieties have many 
distinct components, of various dimensions
(\cite{CerveauDeserti:ptdegre, BisiCalabriMella}). 
The group 
$\mathrm{Aut}(\mathbb{P}^k_\mathbb{C})=\mathrm{PGL}(k+1,\mathbb{C})$ 
acts by left translations, by right translations, 
and by conjugacy on 
$\mathrm{Bir}_d(\mathbb{P}^k_\mathbb{C})$. 
Since this group is connected, these actions preserve each connected component. 

A birational map 
$f\colon\mathbb{P}^k_\mathbb{C}\dashrightarrow\mathbb{P}^k_\mathbb{C}$
is \textsl{regularizable} if there there 
exist a smooth projective variety~$V$ 
and a birational map 
$g\colon V\dashrightarrow\mathbb{P}^k_\mathbb{C}$
such that $g^{-1}\circ f\circ g$ is an 
automorphism of $V$.
To any element~$f$ of $\mathrm{Bir}(\mathbb{P}^k_{\mathbb{C}})$ we 
associate the set $\mathrm{Reg}(f)$ defined by
\[ 
\mathrm{Reg}(f):=\big\{ A\in \mathrm{Aut}(\mathbb{P}^k_{\mathbb{C}})\; \vert \; A\circ f \; {\text{ is regularizable}} \big\}.
\]

On the one hand Dolgachev asked whether there exists a 
birational self map of $\mathbb{P}^k_\mathbb{C}$ 
of degree $>1$ such that 
$\mathrm{Reg}(f)=\mathrm{Aut}(\mathbb{P}^k_\mathbb{C})$.
In \cite{SJJ} we give a negative answer 
to this question. More precisely we prove 

\begin{thm}[\cite{SJJ}]
Let $f$ be a birational self map of 
$\mathbb{P}^k_\mathbb{C}$ of degree $d\geq 2$.

The set of automorphisms~$A$ of
$\mathbb{P}^k_\mathbb{C}$ such that 
$\deg\big((A\circ f)^n\big)\not=\big(\deg(A\circ f)\big)^n$
for some $n>0$ is a countable union of 
proper Zariski closed subsets
of $\mathrm{PGL}(k+1,\mathbb{C})$.

In particular there exists an automorphism~$A$ of 
$\mathbb{P}^k_\mathbb{C}$ such that 
$A\circ f$ is not regularizable.
\end{thm}

On the other hand Bedford asked: does there exist a birational 
map $f$ of $\mathbb{P}^k_{\mathbb{C}}$ 
such that $\mathrm{Reg}(f)=~\emptyset$~? 
We will focus on the case $k=2$.
According to \cite{BedfordKim, Diller} if $\deg f=2$,
then $\mathrm{Reg}(f)\not=\emptyset$. What about birational maps
of degree $3$ ? Blanc proves that the set
\[
\big\{f\in\mathrm{Bir}_3(\mathbb{P}^2_{\mathbb{C}})\,\vert\, \mathrm{Reg}(f)\not=\emptyset,\,\displaystyle\lim_{n\to +\infty}(\deg (f^n))^{1/n}>1\big\}
\]
is dense in $\mathrm{Bir}_3(\mathbb{P}^2_{\mathbb{C}})$ 
and that its complement 
has codimension $1$ (\emph{see} \cite{Blanc}). 
Blanc also gives a 
positive answer to Bedford question in dimension $2$:
if $\chi\colon\mathbb{P}^2_{\mathbb{C}}\dashrightarrow\mathbb{P}^2_{\mathbb{C}}$ 
is the birational map given by 
\[
\chi\colon(x:y:z)\dashrightarrow\big(xz^5+(yz^2+x^3)^2:yz^5+x^3z^3:z^6\big)
\]
then $\mathrm{Reg}(\chi)=\emptyset$.

\begin{rem}\label{rem:blancdegsup}
Note that $\chi=(x+y^2,y)\circ(x,y+x^3)$ in the affine 
chart 
$z=1$. Indeed Blanc example can be generalized 
as follows:
the birational map given in the affine chart $z=1$ by 
\[
\chi_{n,p}=(x+y^n,y)\circ(x,y+x^p)=\big(x+(y+x^p)^n,y+x^p\big)
\]
satisfies $\mathrm{Reg}(\chi_{n,p})=\emptyset$.
\end{rem}

In this article we prove that there exist 
a birational self map $\psi$ of $\mathbb{P}^2_{\mathbb{C}}$ 
such that $\deg\psi<6$ and $\mathrm{Reg}(\psi)=\emptyset$:

\begin{theoalph}\label{thm:bedford}
If $\psi\colon\mathbb{P}^2_{\mathbb{C}}\dashrightarrow\mathbb{P}^2_{\mathbb{C}}$ 
is the birational map given by 
\[
\psi\colon(x:y:z)\dashrightarrow\big(x^2yz^2-z^5+x^5:x^2(x^2y-z^3):xz(x^2y-z^3)\big),
\]
then $\mathrm{Reg}(\psi)=\emptyset$.
\end{theoalph}

\subsection*{Acknowledgements} I would like to 
thank Serge Cantat
for many interesting discussions.
I am also grateful to the referee who has led 
me to considerably improve the drafting of the 
article.

\section{Proof of Theorem \ref{thm:bedford}}

Let $S$ be a smooth projective surface.
Let $\phi\colon S\dashrightarrow S$ 
be a birational map. This map admits 
a resolution
\[
 \xymatrix{
     & Z\ar[rd]^{\pi_2}\ar[ld]_{\pi_1} & \\
    S\ar@{-->}[rr]_\phi & & S,
  }
\]
where 
$\pi_1\colon Z\to S$ and $\pi_2\colon Z\to S$ are finite 
sequences of blow-ups. 
The resolution is \textsl{minimal} if and only
if no $(-1)$-curve of $Z$ is contracted by both $\pi_1$ 
and $\pi_2$. The \textsl{base-points} $\mathrm{Base}(\phi)$
of $\phi$ are the points blown-up by $\pi_1$, which can be 
points of $S$ or infinitely near points.
The proper base-points 
of $\phi$ are called \textsl{indeterminacy points} of
$\phi$ and form a set denoted $\mathrm{Ind}(\phi)$.
Finally we denote by $\mathrm{Exc}(\phi)$ the set
of curves contracted by $\phi$.

Denote by 
$\mathfrak{b}(\phi)$
the number of base-points of $\phi$; note that 
$\mathfrak{b}(\phi)$ is equal to the difference
of the ranks of $\mathrm{Pic}(Z)$ and 
$\mathrm{Pic}(S)$ and thus equal to 
$\mathfrak{b}(\phi^{-1})$. Let us introduce the 
\textsl{dynamical number of the base-points of 
$\phi$}
\[
\mu(\phi)=\displaystyle\lim_{k\to +\infty}\frac{\mathfrak{b}(\phi^k)}{k}.
\]
Since 
$\mathfrak{b}(\phi\circ\varphi)\leq\mathfrak{b}(\phi)+\mathfrak{b}(\varphi)$ 
for any birational
self map~$\varphi$ of $S$, 
$\mu(\phi)$ is a non-negative real number.
As $\mathfrak{b}(\phi)=\mathfrak{b}(\phi^{-1})$
one gets $\mu(\phi^k)=\vert k\,\mu(\phi)\vert$ 
for any $k\in\mathbb{Z}$. Furthermore if  $Z$ is a smooth projective surface and
$\varphi\colon S\dashrightarrow Z$ a birational map, then for all 
$n\in\mathbb{Z}$
\[
-2\mathfrak{b}(\varphi)+\mathfrak{b}(\phi^n)\leq\mathfrak{b}(\varphi\circ\phi^n\circ\varphi^{-1})\leq 2\mathfrak{b}(\varphi)+\mathfrak{b}(\phi^n);
\]
hence $\mu(\phi)=\mu(\varphi\circ\phi\circ\varphi^{-1})$. 
One can thus state the following 
result:

\begin{lem}[\cite{BlancDeserti:degree}]\label{lem:losange}
The dynamical number of base-points is an 
invariant of conjugation. In parti\-cular if 
$\phi$ is a regularizable birational self map of a smooth projective surface, then 
$\mu(\phi)=0$.
\end{lem}

A base-point $p$ of $\phi$ is a 
\textsl{persistent base-point}
if there exists an integer $N$ such 
that for any $k\geq N$
\begin{itemize}
\item[$\diamond$] $p\in\mathrm{Base}(\phi^k)$
\item[$\diamond$] and $p\not\in\mathrm{Base}(\phi^{-k})$.
\end{itemize}
Let $p$ be a point of $S$ or a 
point infinitely near $S$ such that 
$p\not\in\mathrm{Base}(\phi)$. Consider
a minimal resolution of $\phi$
\[
 \xymatrix{
     & Z\ar[rd]^{\pi_2}\ar[ld]_{\pi_1} & \\
    S\ar@{-->}[rr]_\phi & & S.
  }
\]
Because $p$ is not a base-point of $\phi$
it corresponds via $\pi_1$ to a point 
of $Z$ or infinitely near; using $\pi_2$
we view this point on $S$ again maybe 
infinitely near and denote it 
$\phi^\bullet(p)$.
For instance if $S=\mathbb{P}^2_\mathbb{C}$, 
$p=(1:0:0)$ and $f$ is the birational 
self map of $\mathbb{P}^2_\mathbb{C}$ given by
\[
(z_0:z_1:z_2)\dashrightarrow(z_1z_2+z_0^2:z_0z_2:z_2^2)
\]
the point $f^\bullet(p)$ is not equal
to $p=f(p)$ but is infinitely near to 
it. Note that if $\varphi$ is 
a birational self map of $S$ and 
$p$ is a point of $S$ such that 
$p\not\in\mathrm{Base}(\phi)$, 
$\phi(p)\not\in\mathrm{Base}(\varphi)$, 
then 
$(\varphi\circ\phi)^\bullet(p)=\varphi^\bullet(\phi^\bullet(p))$. One can put an equivalence relation on 
the set of points of $S$ or infinitely 
near $S$: the point $p$ is 
\textsl{equivalent}
to the point $q$ if there exists an 
integer $k$ such that $(\phi^k)^\bullet(p)=q$;
in particular $p\not\in\mathrm{Base}(\phi^k)$
and $q\not\in\mathrm{Base}(\phi^{-k})$. 
Remark that the equivalence class is the 
generalization of set of orbits for 
birational maps.

Let us give the relationship between
the dynamical number of base-points and 
the equi\-valence classes of 
persistent base-points:

\begin{pro}[\cite{BlancDeserti:degree}]
Let $S$ be a smooth projective surface.
Let $\phi$ be a birational self map 
of $S$.

Then $\mu(\phi)$ coincides with the 
number of equivalence classes of 
persistent base-points of~$\phi$. 
In particular $\mu(\phi)$ is an 
integer.
\end{pro}

This interpretation of the dynamical number of base-points allows to prove the  following result that gives a
characterization of regularizable  birational maps:

\begin{thm}[\cite{BlancDeserti:degree}]\label{thm:doublelosange}
Let $\phi$ be a birational self map of a smooth 
projective surface. Then $\phi$ is regularizable if and only if $\mu(\phi)=0$. 
\end{thm}

\smallskip

\subsection{Base-points of $\psi$}

The birational map 
\[
\psi\colon(x:y:z)\dashrightarrow
\big(x^2yz^2-z^5+x^5:x^2(x^2y-z^3):xz(x^2y-z^3)\big)
\]
has only one proper base-point, namely $p_1=(0:1:0)$, 
and all its base-points are in tower that is: the 
nine base-points of $\psi$ that we denote 
$p_1$, $p_2$, $\ldots$, $p_9$ are such that 
$p_i$ is infinitely near to 
$p_{i-1}$ for $2\leq i\leq 9$. We denote by 
$\pi\colon S\to\mathbb{P}^2_\mathbb{C}$ the 
blow up of the $9$ base-points, and still 
write $L_x$ (resp. $\mathcal{C}$) the strict
transform of the line 
$L_x\subset\mathbb{P}^2_\mathbb{C}$ of 
equation $x=0$ (resp. the 
curve of equation $x^2y-z^3=0$) which is 
contracted by $\psi$. We denote by 
$E_i\subset S$ the strict transform of the 
curve obtained by blowing up $p_i$. The 
configuration of the curves $E_1$, $E_2$, 
$\ldots$, $E_9$, $L_x$ and $\mathcal{C}$
is

\bigskip

\begin{center}
\begin{figure}[H]
\setlength{\unitlength}{1cm}
\begin{picture}(6.2,1.5)
\put(-0.8,0){\circle*{0.15}}
\put(-0.8,0.2){$L_x$}
\put(-0.8,0){\line(1,0){0.8}}
\put(0,0){\circle*{0.15}}
\put(0,0.2){$\mathrm{E}_2$}
\put(0,0){\line(1,0){0.8}}
\put(0.8,0){\circle*{0.15}}
\put(0.8,0.2){$\mathrm{E}_3$}
\put(0.8,0){\line(-1,-1){0.6}}
\put(0.2,-0.6){\circle*{0.15}}
\put(0,-1.1){$\mathrm{E}_1$}
\put(0.8,0){\line(1,0){0.8}}
\put(1.6,0){\circle*{0.15}}
\put(1.6,0.2){$\mathrm{E}_4$}
\put(1.6,0){\line(1,0){0.8}}
\put(2.4,0){\circle*{0.15}}
\put(2.4,0.2){$\mathrm{E}_5$}
\put(2.4,0){\line(1,0){0.8}}
\put(3.2,0){\circle*{0.15}}
\put(3.2,0.2){$\mathrm{E}_6$}
\put(3.2,0){\line(1,0){0.8}}
\put(4,0){\circle*{0.15}}
\put(4,0.2){$\mathrm{E}_7$}
\put(4,0){\line(1,0){0.8}}
\put(4.8,0){\circle*{0.15}}
\put(4.8,0.2){$\mathrm{E}_8$}
\put(4.8,0){\line(1,-1){0.6}}
\put(5.4,-0.6){\circle*{0.15}}
\put(5.4,-1.1){$\mathcal{C}$}
\put(4.8,0){\line(1,0){0.8}}
\put(5.6,0){\circle*{0.15}}
\put(5.6,0.2){$\mathrm{E}_9$}
\end{picture}
\bigskip
\bigskip
\bigskip
\bigskip
\caption{}
\end{figure}
\end{center}

Two curves are connected by an edge if their
intersection is positive. Let us write
$\psi_A=A\circ\psi$ where $A$ is an 
automorphism of $\mathbb{P}^2_\mathbb{C}$.
Because $\pi$ is the blow-up of the 
base-points of $\psi$, which are also the
base-points of $\psi_A$, the map 
$\eta=\psi_A\circ\pi$ is a birational 
morphism $S\to\mathbb{P}^2_\mathbb{C}$ 
which is the blow-up of the 
base-points of $\psi_A^{-1}$. In 
fact 
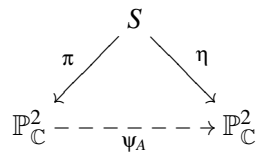
\[
\xymatrix{& S\ar[rd]^{\eta}\ar[ld]_{\pi}&\\
\mathbb{P}^2_\mathbb{C}\ar@{-->}[rr]_{\psi_A}&&\mathbb{P}^2_\mathbb{C}
}
\]
is the minimal resolution of 
$\psi_A$.

The morphism $\eta$ contracts $L_x$ and 
$\mathcal{C}$ as well as the union 
of eight other irreducible curves 
which are among the curves 
$E_1$, $E_2$, $\ldots$, $E_9$. 
The configuration of Figure 1 
shows that $\eta$ contracts the 
curves $L_x$, $E_2$, $E_3$, $E_4$, 
$E_5$, $E_6$, $E_7$, $E_8$, 
$\mathcal{C}$ following this order.

We can see 
$\eta\colon S\to\mathbb{P}^2_\mathbb{C}$
as a sequence of nine blow-ups in the same 
way as we did for $\pi$. We denote 
by $q_1$, $q_2$, $\ldots$, $q_9$ the 
base-points of $\psi_A^{-1}$ (or 
equivalently the points blown up by 
$\eta$) so that
$q_1\in\mathbb{P}^2_\mathbb{C}$
and $q_i$ is infinitely near to 
$q_{i-1}$ for $2\leq i\leq 9$. We 
denote by $D\subset\mathbb{P}^2_\mathbb{C}$
(resp. 
$\mathcal{C}'\subset\mathbb{P}^2_\mathbb{C}$) 
the line contracted by $\psi_A^{-1}$
which is the image by $A$ of the 
line $y=0$ (resp. of 
the conic $z^2-xy=0$). We denote 
by $F_i\subset S$ the strict transform
of the curve obtained by blowing up 
$q_i$. Because of the order of the 
curves contracted by $\eta$ we get 
equalities between $L_x$, $\mathcal{C}$, 
$E_1$, $E_2$, $\ldots$, $E_9$ and 
$D$, $\mathcal{C}'$, $F_1$, $F_2$, 
$\ldots$, $F_9$ as follows

\begin{center}
\begin{small}
\begin{figure}[H]
\setlength{\unitlength}{1.6cm}
\begin{picture}(6.2,1.5)
\put(0,0){\circle*{0.08}}
\put(-0.4,-0.4){$\mathrm{E}_9=L_y$}
\put(0,0){\line(1,0){0.8}}
\put(0.8,0){\circle*{0.08}}
\put(0.6,0.2){$\mathrm{E}_8=\mathrm{F}_2$}
\put(0.8,0){\line(-1,1){0.6}}
\put(0.2,0.6){\circle*{0.08}}
\put(-0.2,0.8){$\mathcal{C}=\mathrm{F}_1$}
\put(0.8,0){\line(1,0){0.8}}
\put(1.6,0){\circle*{0.08}}
\put(1.4,-0.4){$\mathrm{E}_7=\mathrm{F}_3$}
\put(1.6,0){\line(1,0){0.8}}
\put(2.4,0){\circle*{0.08}}
\put(2.2,0.2){$\mathrm{E}_6=\mathrm{F}_4$}
\put(2.4,0){\line(1,0){0.8}}
\put(3.2,0){\circle*{0.08}}
\put(3,-0.4){$\mathrm{E}_5=\mathrm{F}_5$}
\put(3.2,0){\line(1,0){0.8}}
\put(4,0){\circle*{0.08}}
\put(3.8,0.2){$\mathrm{E}_4=\mathrm{F}_6$}
\put(4,0){\line(1,0){0.8}}
\put(4.8,0){\circle*{0.08}}
\put(4.6,-0.4){$\mathrm{E}_3=\mathrm{F}_7$}
\put(4.8,0){\line(1,1){0.6}}
\put(5.4,0.6){\circle*{0.08}}
\put(5,0.8){$\mathrm{E}_1=\mathcal{C}'$}
\put(4.8,0){\line(1,0){0.8}}
\put(5.6,0){\circle*{0.08}}
\put(5.4,0.2){$\mathrm{E}_2=\mathrm{F}_8$}
\put(5.6,0){\line(1,0){0.8}}
\put(6.4,0){\circle*{0.08}}
\put(6.4,-0.4){$L_x=\mathrm{F}_9$}
\end{picture}
\bigskip
\bigskip
\bigskip
\bigskip
\caption{}
\end{figure}
\end{small}
\end{center}

In particular we see that the 
configuration of the points $q_1$, $q_2$, 
$\ldots$, $q_9$ is not the same as that 
of the points $p_1$, $p_2$, $\ldots$, 
$p_9$. Saying that a point $m$ is 
proximate to a point $m'$ if $m$ is 
infinitely near to $m'$ and that it 
belongs to the strict transform of the
curve obtained by blowing up $m'$
the configurations of the points $p_i$
and $q_i$ are 

\bigskip

\begin{center}
$\xymatrix{p_1 & \ar[l] p_2 &\ar[l] \ar@/^1pc/[ll] p_3 & \ar[l] p_4& \ar[l] p_5& \ar[l] p_6& \ar[l] p_7& \ar[l] p_8& \ar[l] p_9}$
\end{center}

\bigskip
\bigskip

\begin{center}
$\xymatrix{q_1 & \ar[l] q_2 &\ar[l] q_3 & \ar[l] q_4& \ar[l] q_5& \ar[l] q_6& \ar[l] q_7& \ar[l] q_8& \ar[l] q_9}$
\end{center}

\bigskip

\begin{center}
\begin{small}
\textsc{Figure} 3
\end{small}
\end{center}

\bigskip

We will prove that for any integer
$i>0$ the point $p_3$ belongs to 
$\mathrm{Base}(\psi_A^i)$ and 
does not belong to $\mathrm{Base}(\psi_A^{-i})$.
It implies that $\mu(\psi_A)>0$
and that $\psi_A$ is not regularizable.

Denote by $k$ the lowest positive
integer such that $p_1$ belongs to
$\mathrm{Base}(\psi_A^{-k})$. If no
such integer exists we write 
$k=\infty$. For any $1\leq i< k$
the point $p_1$ does not belong 
to $\mathrm{Base}(\psi_A^{-i})$ so 
$\psi_A$ and $\psi_A^{-1}$ have no 
common base-point. As a consequence
the set of base-points of the map
$\psi_A^{i+1}=\psi_A\circ\psi_A^i$ 
is the union of the base-points of 
$\psi_A^i$ and of the points 
$(\psi_A^{-i})^\bullet(p_j)$ for 
$1\leq j\leq 9$. Since the map 
$\psi_A^{-i}$ is defined at $p_1$ the 
point $(\psi_A^{-i})^\bullet(p_j)$ 
is proximate to the point 
$(\psi_A^{-i})^\bullet(p_k)$ if and 
only if $p_j$ is proximate to 
$p_k$. Proceeding by induction on $i$
we get the following assertions:
\begin{itemize}
\item[$\diamond$] for any $1\leq i\leq k$ integer
$\mathrm{Base}(\psi_A^i)=\{(\psi_A^{-m})^\bullet(p_j)\,\vert\,1\leq j\leq 9,\, 0\leq m\leq i-1\}$;

\item[$\diamond$] for any $0\leq -\ell\leq k$ 
the configuration of the points 
$\{(\psi_A^\ell)^\bullet(p_j)\,\vert\,1\leq j\leq 9\}$ 
is given by 

\bigskip

\begin{center}
\scalebox{0.6}{$\xymatrix{(\psi_A^\ell)^\bullet(p_1) & \ar[l](\psi_A^\ell)^\bullet(p_2) &\ar[l] \ar@/^1pc/[ll] (\psi_A^\ell)^\bullet(p_3) & \ar[l] (\psi_A^\ell)^\bullet(p_4)& \ar[l] (\psi_A^\ell)^\bullet(p_5)& \ar[l] (\psi_A^\ell)^\bullet(p_6)& \ar[l] (\psi_A^\ell)^\bullet(p_7)& \ar[l] (\psi_A^\ell)^\bullet(p_8)& \ar[l] (\psi_A^\ell)^\bullet(p_9)}$}
\end{center}

\bigskip

\end{itemize}

Hence the point $p_3$ belongs to 
$\mathrm{Base}(\psi_A^i)$ for 
any $1\leq i\leq k$. 

If $k=\infty$, then $p_3$ belongs 
to $\mathrm{Base}(\psi_A^i)$ for 
any $i>0$ and by definition of 
$k$ the point $p_1$ does not 
belong to $\mathrm{Base}(\psi_A^{-i})$
for any $i>0$, and so neither $p_3$.
We can thus assume that $k$ is 
a positive integer.

Assume that $q_1$ belongs to 
$\mathrm{Base}(\psi_A^i)$ for some 
$1\leq i\leq k-1$. Then $q_1$ is 
equal to $(\psi_A^{-m})^\bullet(p_j)$
for some $0\leq m\leq k-2$ and 
$1\leq j\leq 9$. This implies
that $p_j$ belongs to 
$\mathrm{Base}(\psi_A^{m+1})$ which 
is impossible because $m+1\leq k-1$. 
Hence $q_1$ does not belong to 
$\mathrm{Base}(\psi_A^i)$ for 
any $1\leq i\leq k-1$.

We thus see that $\psi_A^{-1}$ has no 
common base-point with $\psi_A^i$ 
for $1\leq i\leq k-1$. In particular
if $B$ denotes 
$\mathrm{Base}(\psi_A^{-1})\cap \mathrm{Base}(\psi_A^k)$,
then
\[
B=\{(\psi_A^{-(k-1)})^\bullet(p_j)\,\vert\,1\leq j\leq 9\}\cap\{q_j\,\vert\,1\leq j\leq 9\}.
\]

Let us remark that $p_1$ belongs to
$\mathrm{Base}(\psi_A^{-k})$ and 
$p_1$ does not belong to $\mathrm{Base}(\psi_A^{-(k-1)})$;
as a result $(\psi_A^{-(k-1)})^\bullet(p_1)$, 
which is a base-point of $\psi_A^k$,
is also a base-point of $\psi_A^{-1}$.
The set $B$ is thus not empty.

The configurations of the two sets 
of points 
$\{(\psi_A^{-(k-1)})^\bullet(p_j)\,\vert\,1\leq j\leq 9\}$
and $\{q_j\,\vert\,1\leq j\leq 9\}$
imply that 
$q_1=(\psi_A^{-(k-1)})^\bullet(p_1)$.

Moreover either $B=\{q_1\}$, or 
$B=\{q_1,\,q_2\}$. Indeed 
$(\psi_A^{-(k-1)})^\bullet(p_3)$ is 
proximate to
$(\psi_A^{-(k-1)})^\bullet(p_2)$
and
$(\psi_A^{-(k-1)})^\bullet(p_1)$ 
whereas
$q_3$ is proximate to $q_2$ but not 
to $q_1$.

The point 
$(\psi_A^{-(k-1)})^\bullet(p_3)$
is thus a point infinitely near to 
$q_1$ in the second neighborhood 
which is maybe infinitely near to $q_2$ 
but not equal to $q_3$. Recalling 
that $\eta$ is the blow up of 
$q_1$, $q_2$, $\ldots$, $q_9$ the 
point 
$(\eta^{-1}\circ\psi_A^{-(k-1)})^\bullet(p_3)$ 
corresponds to a point that belongs, 
as a proper or infinitely near point, 
to one of the curves $F_1$, 
$F_2\subset S$.
So 
$(\pi\circ\eta^{-1}\circ\psi_A^{-(k-1)})^\bullet(p_3)$
is a point infinitely near to $p_3$.
For any $1\leq i\leq k$ the point 
$p_3$ does not belong to 
$\mathrm{Base}(\psi_A^{-i})$; 
therefore there is no base-point of 
$\psi_A^{-i}$ which is infinitely
near to $p_3$. As a result
$(\psi_A^{-k})^\bullet(p_3)$ does not 
belong to $\mathrm{Base}(\psi_A^{-i})$
and $p_3$ does not belong to 
$\mathrm{Base}(\psi_1^{-(k+i)})$. 
Moreover $(\psi_A^{-(k+i)})^\bullet(p_3)$
is infinitely near to
$(\psi_A^{-i})^\bullet(p_3)$. 
Choosing $i=k$ we see that 
$(\psi_A^{-2k})^\bullet(p_3)$ is
infinitely near to
$(\psi_A^{-k})^\bullet(p_3)$ which is 
infinitely near to $p_3$. Continuing
like this we get
\[
\forall\,i\geq 1\qquad p_3\not\in\mathrm{Base}(\psi_A^{-i}).
\]
To get the result it remains to show
that $p_3$ belongs to 
$\mathrm{Base}(\psi_A^i)$ for any 
$i\geq 1$. Reversing the order of $\psi_A$ and 
$\psi_A^{-1}$ we prove as previously that
\[
\forall\,i\geq 1\qquad q_3\not\in\mathrm{Base}(\psi_A^i).
\]

Let us now see that 
\[
\big(\forall\,i\geq 1 \quad q_3\not\in\mathrm{Base}(\psi_A^i)\big)\Rightarrow\big(\forall\,i\geq 1 \quad p_3\in\mathrm{Base}(\psi_A^i)\big).
\]
For $i=1$ it is obvious. Assume $i>1$; let us decompose 
\begin{itemize}
\item[$\diamond$] $\psi_A^i$ into $\psi_A^{i-1}\circ\psi_A$,
\item[$\diamond$] $\pi\colon S\to\mathbb{P}^2_\mathbb{C}$ into 
$\pi_{12}\circ\pi_{39}$ where $\pi_{12}\colon Y\to\mathbb{P}^2_\mathbb{C}$
is the blow up of $p_1$, $p_2$ and $\pi_{39}\colon S\to~Y$
is the blow up of $p_3$, $p_4$, $\ldots$, $p_9$,
\item[$\diamond$] $\eta\colon S\to\mathbb{P}^2_\mathbb{C}$ into 
$\eta_{12}\circ\eta_{39}$ where $\eta_{12}\colon Z\to\mathbb{P}^2_\mathbb{C}$
is the blow up of $q_1$, $q_2$ and $\eta_{39}\colon S\to~Z$
is the blow up of $q_3$, $q_4$, $\ldots$, $q_9$. 
\end{itemize}
Note that $\eta_{39}$ contracts $F_9$, $F_8$, $\ldots$, $F_3$ 
onto the point $Z\ni q_3\not\in\mathrm{Base}(\psi_A^{i-1}\circ\eta_{12})$. 
Consider the system of conics of $\mathbb{P}^2_\mathbb{C}$
passing through $p_1$, $p_2$ and $p_3$. Denote by $\Lambda$ 
its lift on $Y$; it is a system of smooth curves 
passing through $q_3$ with movable tangents and 
$\dim\Lambda=2$. The strict transform of $\Lambda$ on $S$ is 
a system of curves intersecting $E_3$ at a general movable
point. The map $\eta_{39}$ contracts the curves $L_x$, $E_2$, 
$E_3$, $E_4$, $E_5$, $E_6$, $E_7$. As the curve $E_3$ is 
contracted and is not the last one, the image of the 
system by $\eta_{39}$ passes through $q_3$ with a fixed 
tangent corresponding to the point $q_4$. Since 
$q_3\not\in\mathrm{Base}(\psi_A^{i-1}\circ\eta_{12})$ the
image of $\Lambda\subset Y$ by 
$\psi_A^{i-1}\circ\eta\circ(\pi_{39})^{-1}$ has a fixed 
tangent at the point $(\psi_A^{i-1}\circ\eta_{12})(q_3)$.
As a consequence $p_3$ belongs to 
$\mathrm{Base}(\psi_A^{i-1}\circ\eta\circ(\pi_{39})^{-1})$
and thus to 
$\mathrm{Base}(\underbrace{\psi_A^{i-1}\circ\eta\circ(\pi_{39})^{-1}\circ(\pi_{12})^{-1}}_{\psi_A^i})$.

\bibliographystyle{plain}

\bibliography{biblio}

\end{document}